\documentclass[11pt]{amsart}

\usepackage{amsmath, amssymb}
\usepackage{mathrsfs}
\usepackage{amsfonts}
\usepackage[all]{xy}
\usepackage{epsfig}

\numberwithin{equation}{section}

\DeclareMathOperator{\Aut}{Aut}
\DeclareMathOperator{\End}{End}
\DeclareMathOperator{\Hom}{Hom}

\DeclareMathOperator{\Ker}{Ker}
\DeclareMathOperator{\Mod}{mod-}
\DeclareMathOperator{\HN}{HN}

\def\lam{\lambda}
\def\endproof{\hfill{$\square$} \medskip}
\def\C{\mathbb{C}}
\def\R{\mathbb{R}}
\def\Z{\mathbb{Z}}
\def\alp{\alpha}
\def\H{\mathcal{H}}
\def\lto{\longrightarrow}
\def\iso{\xrightarrow{\sim}}

\def\sli{\mathcal{P}}
\def\F{\mathcal{F}}
\def\D{\mathcal{D}}
\def\K{\mathbb{K}}
\def\Stab{\mathrm{Stab}}
\def\Sph{\mathrm{Sph}}
\def\bH{\mathbb{H}}
\def\Q{\mathcal{Q}}
\def\eps{\epsilon}
\def\Db{\mathrm{D}^b}

\newtheorem{thm}{Theorem}[section]
\newtheorem{prop}[thm]{Proposition}

\newtheorem{rem}[thm]{Remark}
\newtheorem{lem}[thm]{Lemma}

\theoremstyle{definition}
\newtheorem{defi}[thm]{Definition}

\renewcommand{\Im}{\text{Im\,}}

\newcommand{\norm}[1]{\lVert #1 \rVert}
\newcommand{\abv}[1]{\lvert #1 \rvert}

\author{Akishi Ikeda}
\address{Kavli Institute for the Physics and Mathematics of the Universe (WPI), UTIAS, 
The University of Tokyo, Kashiwa, Chiba 277-8583, Japan}
\email{akishi.ikeda@ipmu.jp}

\begin{document}

\title[Mass growth of objects and categorical entropy]
{Mass growth of objects and categorical entropy}

\begin{abstract}
In the pioneer work by Dimitrov-Haiden-Katzarkov-Kontsevich, they introduced
various categorical analogies from  classical theory of 
dynamical systems. In particular, they 
defined the entropy of an endofunctor 
on a triangulated category 
with a split generator. In the connection between categorical theory 
and classical theory, 
a stability condition on a 
triangulated category plays the role of a measured foliation 
so that one can measure the ``volume'' of objects, called the mass, 
via the stability condition. 
The aim of this paper is to establish fundamental properties of 
the growth rate of mass of objects 
under the mapping by the endofunctor and to clarify 
the relationship between  the 
entropy and that. We also show that they coincide 
under a certain condition. 
\end{abstract}

\maketitle

\section{Introduction}
In the pioneer work \cite{DHKK}, Dimitrov-Haiden-Katzarkov-Kontsevich introduced
various categorical analogies of classical theory from 
dynamical systems. In particular, they 
defined the entropy 
of an endofunctor on a triangulated category 
with a split generator. 
One of their motivations comes from the connection 
between theory of stability conditions on triangulated categories 
and Teichm\"uller theory of surfaces \cite{GMN,BrSm}. 
In this connection, a stability condition on a 
triangulated category corresponds 
to a measured foliation (a quadratic differential) on a surface, 
and the mass of stable objects  
corresponds to the length of  geodesics. 
Thus the mass of objects plays the role of ``volume'' in some sense. 
In the work \cite{DHKK}, they also suggested that there is a connection 
between the growth rate of mass of objects under the mapping 
by an endofunctor and the entropy of that. 
In this paper, we establish fundamental properties of 
the mass growth and clarify the relationship 
between the entropy and that. We also show that they 
coincide under a certain condition.
The result in this paper is motivated by 
the famous classical work ``Volume growth and entropy'' 
by Yomdin \cite{Yom} on classical dynamical systems.

\subsection{Fundamental properties of mass growth}
First we introduce the mass growth with respect to endofunctors. 
Let $\D$ be a triangulated category and $K(\D)$ be its Grothendieck group. 
A stability condition $\sigma=(Z,\sli)$ on $\D$ \cite{Br1} 
is a pair of a linear map $Z \colon K(\D) \to \C$
and a family of full subcategories $\sli(\phi) \subset \D$ for $\phi \in \R$ 
satisfying some axioms (see Definition \ref{def:stab}). 
A nonzero object in $\sli(\phi)$ is called a semistable object of phase $\phi$. 
One of the axioms implies that
any nonzero object $E \in \D$ can be decomposed  
into semistable objects with decreasing phases, i.e. 
there is a sequence of exact triangles, called a Harder-Narasimhan filtration,
\begin{equation*}
0 =
\xymatrix @C=4.5mm{ 
 E_0 \ar[rr]   &&  E_1 \ar[dl] \ar[rr] && E_2 \ar[dl] 
 \ar[r] & \dots  \ar[r] & E_{m-1} \ar[rr] && E_m \ar[dl] \\
& A_1 \ar@{-->}[ul] && A_2 \ar@{-->}[ul] &&&& A_m \ar@{-->}[ul] 
}
= E
\end{equation*}
with $A_i \in \sli(\phi_i)$ and $\phi_1>\phi_2 >\cdots>\phi_m$. 
Through the Harder-Narasimhan filtration, 
the mass of $E$ with a parameter $t \in \R$ 
(see Definition \ref{defi:mass}) is defined by
\[
m_{\sigma,t}(E):=\sum_{i=1}^m |Z(A_i)|e^{\phi_i t}. 
\] 
Thus a given stability condition defines the ``volume'' of objects in some sense. 
Actually in the connection between spaces of stability conditions and 
Teichm\"uller spaces, the mass of stable objects gives 
the length of corresponding geodesics \cite{BrSm,GMN,HKK,Ike}. 
For an endofunctor $F \colon \D \to \D$, we want to consider 
the growth rate of mass of objects under the mapping by $F$. 
Therefore we introduce the following quantity.
The mass growth with respect to $F$ is 
the function $h_{\sigma,t}(F) \colon \R \to [-\infty,\infty]$
defined by
\[
h_{\sigma,t}(F):=\sup_{E \in \D} \left\{ \,\limsup_{n \to \infty}\frac{1}{n}\,
\log (m_{\sigma,t}(F^n E))  
\, \right\}.
\]
(As conventions, set $m_{\sigma,t}(0)=0$ and $\log 0 = -\infty$.)
Fundamental properties of $h_{\sigma,t}(F)$ are stated 
as the main result of this paper. 
We also recall the space of stability conditions
to consider the behavior of $h_{\sigma,t}(F)$ under the 
deformation of $\sigma$. 
In \cite{Br1}, it was shown that the set of stability conditions 
$\Stab(\D)$ has a natural topology  
and in addition, $\Stab(\D)$ becomes a complex manifold. 

Next we recall the entropy  of endofunctors from \cite{DHKK}. 
Let $\D$ be a triangulated category with a split-generator and 
$F \colon \D \to \D$ be an endofunctor. 
In \cite{DHKK}, they introduced the function 
$h_t(F) \colon \R \to [-\infty,\infty)$, called the entropy of $F$ 
(see Definition \ref{def:entropy}), 
and showed various fundamental properties of $h_t(F)$.  
In addition, they asked the relationship between 
the entropy $h_t(F)$ and the mass growth $h_{\sigma,t}(F)$ 
(see \cite[Section 4.5]{DHKK}). 
Our result is the following. 
\begin{thm}[Theorem \ref{thm:well} and Proposition \ref{prop:connected}]
Let $\D$ be a triangulated category, $F \colon \D \to \D$ be an endofunctor and 
$\sigma$ be a stability condition on $\D$.
Assume that $\D$ has 
a split-generator $G$. Then the mass growth
$h_{\sigma,t}(F)$ satisfies the followings. 
\begin{itemize}
\item[(1)] If a stability condition $\tau$ lies in the same connected 
component as $\sigma$ in the space of stability conditions $\Stab(\D)$, then 
\[
h_{\sigma,t}(F)=h_{\tau,t}(F).
\]
\item[(2)] The mass growth of the generator $G$ determines
$h_{\sigma,t}(F)$, i.e.
\[
h_{\sigma,t}(F)=\limsup_{n \to \infty}\frac{1}{n}\,
\log (m_{\sigma,t}(F^n G)).  
\]
\item[(3)] 
An inequality 
\[
h_{\sigma,t}(F) \le h_t(F)<\infty
\]
holds. 
\end{itemize}
\end{thm}
In the case $t=0$, this result was stated in \cite[Section 4.5]{DHKK} 
by using the triangle inequality for mass 
(see Proposition \ref{prop:mass_triangle}). 
However, the triangle inequality for 
mass is non-trivial even if $t=0$
and actually the most technical part in this paper. 
Therefore we give  
a detailed proof of it with a parameter $t$ in Section \ref{sec:tri}.

\subsection{Lower bound by the spectral radius}
We consider the lower bound of the mass growth
when $t=0$. 
Since $F \colon \D \to \D$ preserves exact triangles, $F$ induces 
a linear map $[F] \colon K(\D) \to K(\D)$. 
The spectral radius of $[F]$ is defined by 
\[
\rho([F]):=\max \{|\lam|\,|\, \lam \text{ is an eigenvalue of }[F]  \,   \}. 
\]
\begin{thm}[Proposition \ref{prop:spectral}]
\label{thm:bound}
In the case $t=0$, we have an inequality
\[
\log \rho([F]) \le h_{\sigma,0}(F) \le h_0(F)
\]
for any stability condition $\sigma \in \Stab(\D)$.
\end{thm}
As known results, if $\D$ is saturated, then 
it was shown in \cite[Theorem 2.9]{DHKK} that  
for a linear map $\mathrm{HH}_*(F) 
\colon \mathrm{HH}_*(\D) \to \mathrm{HH}_*(\D)$
induced on the Hochschild homology of $\D$, the inequality 
$\log \rho (\mathrm{HH}_*(F)) \le h_0(F)$ holds under some  
condition for eigenvalues of  $\mathrm{HH}_*(F)$. 
They also conjectured that the inequality 
holds without that condition.  
Our result Theorem \ref{thm:bound} 
holds without any conditions for $[F]$,   
however we use the existence of stability conditions on $\D$. 
For many examples in \cite{DHKK,Kik,KT}, 
it was shown that the equality $\log \rho ([F]) = h_0(F)$ holds.
Kikuta-Takahashi gave a certain conjecture on the equality  
in \cite[Conjecture 5.3]{KT}.

\subsection{Equality between mass growth and entropy}
 The remaining important question is to ask when the equality 
 $h_{\sigma,t}(F)=h_t(F)$ holds. 
 In the following, we give a sufficient condition for the equality. 
 For a stability condition $\sigma=(Z,\sli)$, we can associate an abelian category, 
 called the heart of $\sli$, 
as the extension-closed subcategory 
 generated by objects in $\sli(\phi)$ for $\phi \in (0,1]$. Denote it by 
 $\sli((0,1])$. 
 A stability condition $\sigma=(Z,\sli)$ is called algebraic if 
 the heart $\sli((0,1])$ is a finite length abelian category with 
 finitely many simple objects 
 (see Definition \ref{defi:algebraic_heart} and Definition \ref{defi:algebraic_stab}). 
\begin{thm}[Theorem \ref{thm:algebraic}]
\label{thm:intro2}
Let $G \in \D$ be a split-generator and $F \colon \D \to \D$ be 
an endofunctor. 
If a connected component 
$\Stab^{\circ}(\D) \subset \Stab(D)$ contains an algebraic
stability condition, then for any $\sigma \in \Stab^{\circ}(\D)$
we have
\[
h_t(F)=h_{\sigma,t}(F)=\lim_{n \to \infty}\frac{1}{n}\log (m_{\sigma,t}(F^n G)).
\]
\end{thm}
 Note that in the above theorem, a stability condition $\sigma$ is not 
 necessarily an algebraic stability condition. 
 
 We see a typical example which satisfies 
 the condition in Theorem \ref{thm:intro2} from Section \ref{sec:dga}. 
 Let $A=\oplus_{k}A^k $ be a dg-algebra such 
 that $H^0(A) $ is a finite dimensional algebra and $H^k(A)=0$ for $k>0$. 
 Denote by $\D_{fd}(A)$ the derived category of dg-modules over $A$   
 with finite dimensional total cohomology, i.e. $\sum_k \dim H^k(M)<\infty $. 
 Then there is a bounded t-structure whose heart is isomorphic to 
 the abelian category of finite dimensional modules over $H^0(A)$. 
As a result, we can construct algebraic stability conditions on $\D_{fd}(A)$. 
Thus in the context of representation theory, Theorem \ref{thm:intro2} 
works well. As an application, we compute the entropy of spherical twists 
in Section \ref{sec:spherical}. 

On the other hand, for derived categories coming from algebraic geometry, 
we cannot find algebraic hearts in general. 
Only in special cases, for exmaple in the case that the derived category has 
a full strong exceptional collection, 
the work by Bondal \cite{Bo} enables us to find algebraic hearts.   
It is an important problem to answer whether the equality 
$h_{\sigma,t}(F)=h_t(F)$ holds without the existence of algebraic 
stability conditions. 
 
 \subsection{Categorical theory versus classical theory}
We compere our result with the famous classical result
``Volume growth and entropy''  by Yomdin \cite{Yom}.  
Let $M$ be a compact smooth manifold and 
$f \colon M \to M$ be a smooth map. The map $f$ induces 
a linear map $f_* \colon H_*(M;\R) \to H_*(M;\R)$ on the 
homology group $H_*(M;\R)$. For the map $f$, we can define 
the topological entropy $h_{top}(f)$ \cite{AKM} and the inequality 
$\log \rho (f_*) \le h_{top}(f)$ was conjectured in \cite{Shu}. 
In \cite{Yom}, Yomdin introduced the 
 volume growth $v(f)$ by using a Riemannian metric on $M$ and showed that
 \[
 \log \rho (f_*) \le v(f)\le h_{top}(f).
 \] 
 Our result Theorem \ref{thm:bound} looks like the categorical analogy of 
 this classical result. On the other hand, the difference between categorical theory 
 and classical theory is that the categorical entropy $h_t(F)$ and 
 the mass growth $h_{\sigma,t}(F)$ have the parameter $t$ which 
measures the growth rate of degree shifts in a 
 triangulated category. This point is an essentially new feature of 
 categorical theory. 

 \subsection*{Acknowledgements}
The author would like to thank Kohei Kikuta,  
Genki Ouchi and Atsushi Takahashi
for valuable discussions and comments. 

This work is supported by World Premier International
Research Center Initiative (WPI initiative), MEXT, Japan, 
JSPS KAKENHI Grant Number JP16K17588 and 
JSPS bilateral Japan-Russia Research Cooperative Program.
This paper was written while the author was visiting Perimeter Institute 
for Theoretical Physics by JSPS Program for Advancing Strategic
International Networks to Accelerate the Circulation of Talented Researchers.
Research at Perimeter Institute 
is supported by the Government of Canada 
through the Department of Innovation, Science and Economic Development Canada and 
by the Province of Ontario through the Ministry of Research, Innovation and Science.

\subsection*{Notations}
We work over a field $\K$. 
All triangulated categories in this paper are  $\K$-linear and 
their Grothendieck groups are free of finite rank, i.e.   
$K(\D) \cong \Z^n$ for some $n$. 
An endofunctor $F \colon \D \to \D$ refers to an exact endofunctor, 
i.e. $F$ preserves all exact triangles and commutes with degree shifts. 
The natural logarithm is extended to 
$\log \colon [0,\infty) \to [-\infty,\infty)$ 
 by setting $\log 0:=-\infty$. 
\section{Preliminaries}
In this section, we prepare basic terminologies 
mainly from \cite{DHKK,Br1}. 
\subsection{Complexity and entropy}
First we recall the notion of 
complexity and entropy from \cite[Section 2]{DHKK}. 

Let $\D$ be a triangulated category. 
A triangulated subcategory is called {\it thick} if 
it is closed under taking direct summands. 
For an object $E \in \D$, we denote by $\langle E \rangle \subset \D$ 
the smallest thick triangulated subcategory 
containing $E$. An object $G \in \D$ is called a {\it split-generator} 
if $\langle G \rangle=\D$.
This implies that for any object $E \in \D$, there 
is some object $E^{\prime} \in \D$ such that we have a sequence of exact triangles
\begin{equation*}
0 =
\xymatrix @C=4.5mm{ 
 A_0 \ar[rr]   &&  A_1 \ar[dl] \ar[rr] && A_2 \ar[dl] 
 \ar[r] & \dots  \ar[r] & A_{k-1} \ar[rr] && A_k \ar[dl] \\
& G[n_1] \ar@{-->}[ul] && G[n_2] \ar@{-->}[ul] &&&& G[n_k] \ar@{-->}[ul] 
}
= E\oplus E^{\prime}
\end{equation*}
with $n_i \in \Z$. We note that the object $E^{\prime}$ and the above 
sequence  
are not unique. 

\begin{defi}[\cite{DHKK}, Definition 2.1]
Let $E_1$ and $E_2$ be objects in $\D$. The {\it complexity of $E_2$ relative to $E_1$} 
is the function $\delta_t(E_1,E_2) \colon \R \to [0,\infty)$ defined by
\[
\delta_t(E_1,E_2):=
\begin{cases}
0 &\text{if }E_2 \cong 0 \\
\inf\left\{ \displaystyle\sum_{i=1}^k e^{n_i t} \,\middle|\,
\begin{xy}
(0,5) *{0}="0", (20,5)*{A_{1}}="1", (30,5)*{\dots}, (40,5)*{A_{k-1}}="k-1", 
(60,5)*{E_2\oplus E_2^{\prime}}="k",
(10,-5)*{E_1[n_{1}]}="n1", (30,-5)*{\dots}, (50,-5)*{E_1[n_{k}]}="nk",
\ar "0"; "1"
\ar "1"; "n1"
\ar@{-->} "n1";"0"
\ar "k-1"; "k" 
\ar "k"; "nk"
\ar@{-->} "nk";"k-1"
\end{xy}
\, \right\} 
&\text{if }E_2 \in \langle E_1\rangle \\
\infty &\text{if }E_2 \notin \langle E_1\rangle.
\end{cases}
\]
\end{defi}

By definition, we have an inequality 
$0<\delta_t(G,E)<\infty$ for a split-generator $G \in \D$ 
and a nonzero object $E \in \D$. 
We recall fundamental inequalities for complexity. 
\begin{prop}[\cite{DHKK}, Proposition 2.3]
\label{prop:DHKK}
For $E_1,E_2,E_3 \in \D$,
\begin{itemize}
\item[(1)] $\delta_t(E_1,E_3) \le \delta_t(E_1,E_2)\delta_t(E_2,E_3)$,
\item[(2)] $\delta_t(E_1,E_2 \oplus E_3) \le \delta_t(E_1,E_2)+\delta_t(E_1,E_3)$,
\item[(3)] $\delta_t(F(E_1),F(E_2))\le\delta_t(E_1,E_2)$ for an 
endofunctor $\D \to \D$. 
\end{itemize}
\end{prop}

Similar to \cite[Proposition 2.3]{DHKK}, 
it is easy to check the following.  
\begin{lem}
\label{lem:triangle}
For objects $D,E_1,E_2,E_3 \in \D$, if there
 is an exact triangle 
$E_1 \to E_2 \to E_3 \to E_1[1]$,  
then 
\[
\delta_t(D,E_2) \le \delta_t(D,E_1)+\delta_t(D,E_3).
\]
\end{lem}

Now we introduce the notion of the entropy of endofunctors. 
The entropy of an endofunctor $F$ 
measures the growth rate of  complexity $\delta_t(G,F^nG)$ 
as $n \to \infty$. 
\begin{defi}[\cite{DHKK}, Definition 2.5]
\label{def:entropy}
Let $\D$ be a triangulated category with a split-generator $G$ and 
let $F \colon \D \to \D$ be an endofunctor. 
The {\it entropy} of $F$
is the function $h_t(F) \colon \R \to [-\infty,\infty)$ defined by
\[
h_t(F):=\lim_{n \to \infty}\frac{1}{n}\log \delta_t(G,F^n G).
\]
\end{defi}
By \cite[Lemma 2.6]{DHKK}, 
it follows that $h_t(F)$ is well-defined and $h_t(F)<\infty$. 

\subsection{Bounded t-structures and the associated cohomology}
\label{sec:t-structure}
\begin{defi}[\cite{BBD}]
A {\it t-structure} on $\D$ is a full subcategory 
$\F \subset \D$ satisfying the following conditions:
\begin{itemize}
\item[(a)] $\F[1] \subset \F$,
\item[(b)] define $\F^{\perp} :=\{F \in \D \vert \Hom(D,F) =0 \text{ for all }D \in \F\,  \}$, 
then for every object $E \in \D$ there is an exact triangle 
$D \to E \to F \to D[1]$ in $\D$ 
with $D \in \F$ and $F \in \F^{\perp}$.
\end{itemize}
In addition, the t-structure $\F \subset \D$ is said to be {\it bounded} 
if $\F$ satisfies the condition
\begin{equation*}
\D = \bigcup_{i,j \in \Z} \F^{\perp}[i] \cap \F[j].
\end{equation*}
\end{defi}

For a t-structure $\F \subset \D$, we define the {\it heart} $\H \subset \D$ by
\begin{equation*}
\H := \F^{\perp}[1] \cap \F.
\end{equation*}
It was proved in \cite{BBD} that $\H$ becomes an abelian category. 
Bridgeland gave the characterization of the heart of a bounded t-structure as follows. 
\begin{lem}[\cite{Br1}, Lemma 3.2]
Let $\H \subset \D$ be a full additive subcategory. 
Then $\H$ is the heart of a bounded t-structure if and only if the 
following conditions hold:
\begin{itemize}
\item[(a)]if $k_1 >k_2 \in \Z$ and $A_i \in \H[k_i]\,(i =1,2)$, 
then $\Hom_{\D}(A_1,A_2) = 0$,
\item[(b)]for $0 \neq E \in \D$, there is a finite sequence of integers
\[
k_1 > k_2 > \cdots > k_m
\]
and a sequence of exact triangles
\[
0 =
\xymatrix @C=4.5mm{ 
 E_0 \ar[rr]   &&  E_1 \ar[dl] \ar[rr] && E_2 \ar[dl] 
 \ar[r] & \dots  \ar[r] & E_{m-1} \ar[rr] && E_m \ar[dl] \\
& A_1 \ar@{-->}[ul] && A_2 \ar@{-->}[ul] &&&& A_m \ar@{-->}[ul] 
}
= E
\]
with $A_i \in \H[k_i]$ for all $i$.
\end{itemize}
\end{lem}
The above filtration in the condition (b) defines the $k$-th cohomology  
$H^k(E) \in \H$ of  the object $E$ by 
\[
H^k(E):=
\begin{cases}
A_i[-k_i] \quad &\text{if }k=-k_i \\
0 &\text{otherwise}.
\end{cases}
\]
This cohomology becomes a cohomological functor from $\D$ to $\H$, 
i.e. if there is an exact triangle $D \to E \to F \to E[1]$, then 
we can obtain a long exact sequence
\[
\cdots \to H^{k-1}(F) \to H^k(D) \to H^k(E) \to H^k(F) \to H^{k+1}(D)\to \cdots
\]
in the abelian category $\H$. 
In the last of this section, we introduce the
special class of bounded t-structures. 
\begin{defi}
\label{defi:algebraic_heart}
We say that the heart of a bounded t-structure is  
 {\it algebraic} if it is a finite length abelian category with 
finitely many isomorphism classes of simple objects. 
\end{defi}
If $\D$ has an algebraic heart $\H$ with simple objects $S_1,\dots,S_n$,
then it is easy to see that the direct sum 
$G:=\oplus_{i=1}^n S_i$ becomes a split-generator of $\D$.

\subsection{Bridgeland stability conditions}
\label{sec:stab}
In \cite{Br1}, Bridgeland introduced the notion of a stability condition 
on a triangulated category as follows.
\begin{defi}
\label{def:stab}
Let $\D$ be a triangulated category and $K(\D)$ be its Grothendieck group. 
A {\it stability condition} $\sigma = (Z, \sli)$ on $\D$ consists of 
a group homomorphism $Z \colon K(\D) \to \C$, called a {\it central charge}, 
and a family of full additive subcategories 
$\sli (\phi) \subset \D$ for $\phi \in \R$ 
satisfying the following conditions:
\begin{itemize}
\item[(a)]
if  $0 \neq E \in \sli(\phi)$, 
then $Z(E) = m(E) \exp(i \pi \phi)$ for some $m(E) \in \R_{>0}$,
\item[(b)]
for all $\phi \in \R$, $\sli(\phi + 1) = \sli(\phi)[1]$, 
\item[(c)]if $\phi_1 > \phi_2$ and $A_i \in \sli(\phi_i)\,(i =1,2)$, 
then $\Hom_{\D}(A_1,A_2) = 0$,
\item[(d)]for $0 \neq E \in \D$, there is a finite sequence of real numbers 
\[
\phi_1 > \phi_2 > \cdots > \phi_m
\]
and a sequence of exact triangles
\[
0 =
\xymatrix @C=4.5mm{ 
 E_0 \ar[rr]   &&  E_1 \ar[dl] \ar[rr] && E_2 \ar[dl] 
 \ar[r] & \dots  \ar[r] & E_{m-1} \ar[rr] && E_m \ar[dl] \\
& A_1 \ar@{-->}[ul] && A_2 \ar@{-->}[ul] &&&& A_m \ar@{-->}[ul] 
}
= E
\]
with $A_i \in \sli(\phi_i)$ for all $i$.
\end{itemize}
\end{defi}

We write $\phi^+_{\sigma}(E):=\phi_1$ and $\phi^-_{\sigma}(E):=\phi_m$.
Nonzero objects in $\sli(\phi)$ are called {\it $\sigma$-semistable
 of phase $\phi$ in $\sigma$}. 
The sequence of exact triangles in (d) is called a {\it Harder-Narasimhan filtration} 
of $E$ with semistable factors $A_1,\dots,A_m$ of phases $\phi_1>\dots>\phi_m$. 

In addition to the above axioms, we always assume that our stability conditions have 
the {\it support property} in \cite{KS}.
Let $\norm{\,\cdot\,}$ be some norm on $K(\D) \otimes \R$. 
A stability condition $\sigma=(Z,\sli)$ 
satisfies the support property if there is a some constant $C >0$ such that 
\begin{equation*}
\frac{\abv{Z(E)}}{\norm{[E]}} >C 
\end{equation*}
for all $\sigma$-semistable objects $E \in \D$. 

For an interval $I \subset \R$, we denote by $\sli(I)$ the extension-closed 
subcategory generated by objects in $\sli(\phi)$ for $\phi \in I$, namely
\[
\sli(I):=\{\,E \in \D \,\vert\, \phi_{\sigma}^{\pm}(E) \in I  \,\} \cup \{0\}. 
\] 
From a stability condition $(Z,\sli)$, we can construct a bonded t-structure 
$\F:=\sli((0,\infty))$ and its heart is given by $\H=\sli((0,1])$.

\subsection{Algebraic stability conditions}
\label{sec:algebraic_stab}
In \cite{Br1}, Bridgeland gave the alternative description of a stability condition on $\D$ 
as a pair of a bounded t-structure and a central charge on its heart. 
By using this description, we construct 
algebraic stability conditions. 

\begin{defi}
Let $\H$ be an abelian category and let $K(\H)$ 
be its Grothendieck group. 
A central charge on $\H$ is 
a group homomorphism $Z \colon K(\H) \to \mathbb{C}$ such that for any nonzero object
$0 \neq E \in \H$, the complex number $Z(E)$ lies in the 
semi-closed upper half-plane 
$\bH := \{\, r e^{i \pi \phi}\ \in \mathbb{C} \,
 \vert \, r \in \mathbb{R}_{>0}, \phi \in (0,1] \,\}$.
\end{defi}
For any nonzero object $E \in \H$, define the {\it phase of $E$} by 
\[
\phi(E):=\frac{1}{\pi} \arg Z(E) \in (0,1]. 
\]
An object $0 \neq E \in \H$ is called {\it $Z$-semistable} if every subobject 
$0 \neq A \subset E$ satisfies $\phi(A) \le \phi(E)$. 
A {\it Harder-Narasimhan filtration} of $0 \neq E \in \H$ is the filtration 
\[
0 =E_0 \subset E_1 \subset \cdots \subset E_{m-1}\subset E_m=E
\]
whose extension factors $F_i:=E_i \slash E_{i-1}$ are $Z$-semistable  
with decreasing phases 
\[
\phi(F_1)>\cdots>\phi(F_m).
\]
A central charge $Z$ is said to have the {\it Harder-Narasimhan property} 
if any nonzero object of $\H$ has a Harder-Narasimhan filtration. 
The following gives the another definition of a stability condition. 
\begin{prop}[\cite{Br1}, Proposition 5.3]
\label{abel_stability}
Giving a stability condition on $\D$ is equivalent 
to giving a heart $\H $ 
of a bounded structure on $\D$ and a central charge on $\H$ 
with the Harder-Narasimhan property. 
\end{prop}
In Proposition \ref{abel_stability}, 
the pair $(Z,\H)$ is constructed from a stability condition $(Z,\sli)$ 
by setting $\H:=\sli((0,1])$.
\begin{defi}
\label{defi:algebraic_stab}
A stability condition $(Z,\sli)$ is called {\it algebraic}
if the corresponding heart $\H=\sli((0,1])$ is algebraic 
(for the definition of algebraic hearts, see Definition \ref{defi:algebraic_heart}). 
\end{defi}
Algebraic stability conditions are constructed 
from algebraic hearts as follows.
Let $\H \subset \D$ be an algebraic 
heart with simple objects $S_1,\dots,S_n$. 
Then the Grothendieck group is given by 
$K(\H) \cong \oplus_{i=1}^n \Z [S_i]$. 
Take  $(z_1,\dots,z_n) \in \bH^n$ and define
the central charge $Z \colon K(\H) \to \C$ by the linear extension of 
$Z(S_i) :=z_i$. Then $Z$ has  
the Harder-Narasimhan property by \cite[Proposition 2.4]{Br1}. 
Thus $(Z,\H)$ becomes a stability condition on $\D$. 

\subsection{Harder-Narasimhan polygons}
\begin{figure}[!ht]
\centering
\includegraphics[scale=1.05]{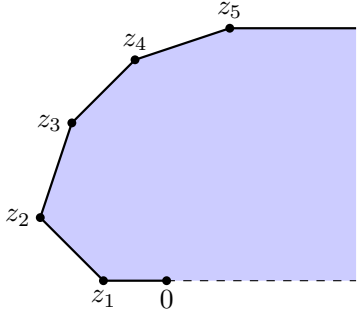}
\caption{Harder-Narasimhan polygon.}
\label{Fig:HN_polygon}
\end{figure}

In this section, we discuss the Harder-Narasimhan polygon 
following \cite{Bay}. 
This plays a key role to show the triangle inequality for mass 
in Section \ref{sec:tri}. 
The following is based on \cite[Section 3]{Bay}. 
\begin{defi}
\label{defi:HN_polygon}
Let $\H$ be an abelian category and $Z$ be a central charge on it. For 
an object $E \in \H$, the {\it Harder-Narasimhan polygon} $\HN^Z(E)$
of $E$ is the convex hull of the subset 
$\{\,Z(A)\in \C\,\vert\, A \subset E\,\} \subset \C$ in the complex plane.
\end{defi}
It is clear from the definition 
that if $F \subset E$, then  $\HN^Z(F) \subset \HN^Z(E)$. 
The Harder-Narasimhan polygon $\HN^Z(E)$ is called {\it polyhedral on the left} 
if it has finitely many extremal points $0=z_0,z_1,\dots,z_k=Z(E)$ such that 
$\HN^Z(E)$ lies to the right of the path $z_0z_1\dots z_k$. 
This implies that the intersection of $\HN^Z(E)$ and the closed half-plane 
to the left of the line through $0$ and $Z(E)$ becomes a polygon 
with vertices $z_0,z_1,\dots,z_k$ (see Figure \ref{Fig:HN_polygon}). 
\begin{prop}[\cite{Bay}, Proposition 3.3]
\label{prop:HN_polygon}
The object $E$ has a Harder-Narasimhan filtration if and only if 
$\HN^Z(E)$  is polyhedral on the left. 
In particular, if the Harder-Narasimhan filtration of $E$ is given by
\[
0=E_0 \subset  E_1 \subset E_2 \subset \cdots E_k=E,
\] 
then extremal points of $\HN^Z(E)$ are given by 
$z_i=Z(E_i)$ for $i=0,1,\dots,k$.
\end{prop}

\subsection{Topology on the space of stability conditions}
In  \cite{Br1}, Bridgeland introduced a natural topology on the space of 
stability conditions and showed that this space becomes a complex manifold. 
In the following, we recall his construction. 
Let $\Stab(\D)$ be the set of stability conditions on 
a triangulated category $\D$ with the support property.  
For stability conditions $\sigma=(Z,\sli)$ and $\tau=(W,\Q)$
in $\Stab(D)$, 
set
\begin{align*}
d(\sli,\Q):= \sup_{0 \neq E \in \D}\left\{\,
\abv{\phi_{\sigma}^-(E) - \phi_{\tau}^-(E)}\,,\, 
\abv{\phi_{\sigma}^+(E) - \phi_{\tau}^+(E)}\,
\right\} \in [0,\infty].
\end{align*}
and
\[
 \norm{Z-W}_{\sigma}:=\sup
 \left\{\frac{\,|Z(E)-W(E)|}{|Z(E)|} \,\middle|\,
 \text{$E$ is $\sigma$-semistable}  \,\right\} \in [0,\infty].
 \]
Define a subset $B_{\eps}(\sigma) \subset \Stab(\D)$ 
by 
\[
B_{\eps}(\sigma):=\{\,\tau=(W,\Q)\in \Stab(\D) \,\vert \,
d(\sli,\Q)<\eps,\, \norm{Z-W}_{\sigma}<\sin(\pi \eps) \, \}
\]
for  $0<\eps<\frac{1}{4}$.

In \cite[Section 6]{Br1}, it was shown that a family of 
subsets
\[
\left\{\,B_{\eps}(\sigma) \subset \Stab(\D) \, \middle\vert\,
\sigma \in \Stab(\D),\,0< \eps<\frac{1}{4} \,\right\}
\]
becomes an open basis of a topology on $\Stab(\D)$. 
In \cite{Br1}, Bridgeland showed a crucial theorem. 
\begin{thm}[\cite{Br1}, Theorem 1.2]
\label{thm:localiso}
The projection map of central charges
\begin{equation*}
\pi \colon \Stab(\D) \longrightarrow \Hom_{\Z}(K(\D),\C), 
\quad (Z,\sli) \mapsto Z
\end{equation*}
 is a local isomorphism of topological spaces. 
 In particular, $\pi$ induces a complex structure on $\Stab(\D)$. 
\end{thm}

\section{Mass growth of objects and categorical entropy}
\subsection{Mass with a parameter and complexity}
In this sectioin, we introduce the mass growth of objects and 
show fundamental properties of it.

\begin{defi}[\cite{DHKK}, Section 4.5]
\label{defi:mass}
Take a stability condition $\sigma=(Z,\sli)$ on $\D$. 
Let $E \in \D$ be  a nonzero object with semistable factors 
$A_1,\dots, A_m$ of phases $\phi_1>\dots>\phi_m$.   
The {\it mass of $E$ with a parameter $t \in \R$} is the function 
$m_{\sigma,t}(E) \colon \R \to \R_{>0}$ defined by
\[
m_{\sigma,t}(E):=\sum_{i=1}^m |Z(A_i)|\,e^{\phi_i t}.
\]
When $t=0$, $m_{\sigma,0}(E)$ 
is called the {\it mass of $E$} and simply written as $m_{\sigma}(E):=m_{\sigma,0}(E)$. 
As a convention, set $m_{\sigma,t}(E):=0$ if $E \cong 0$. 
\end{defi}
In the following, if $\sigma$ is clear in the context, we often drop it 
from the notation and write $m_t(E)$. 
Similar to the growth rate of complexity of a generator 
with respect to endofunctors, 
we consider the growth rate of mass of objects. 
\begin{defi}[\cite{DHKK}, Section 4.5]
Let  $\sigma $ be a stability condition on $\D$ and 
$F \colon \D \to \D$ be an endofunctor. 
The {\it mass growth with respect to $F$} is 
the function $h_{\sigma,t}(F) \colon \R \to [-\infty,\infty]$
defined by
\[
h_{\sigma,t}(F):=\sup_{E \in \D} \left\{ \,\limsup_{n \to \infty}\frac{1}{n}\,
\log (m_{\sigma,t}(F^n E))  
\, \right\}.
\]
\end{defi}
In the rest of this section, we study fundamental properties of $h_{\sigma,t}(F)$. 
The triangle inequality for $m_{\sigma,t}$ plays an important role. 
\begin{prop}
\label{prop:mass_triangle}
For objects $D,E,F \in \D$, 
if there is an exact triangle $D \to E \to F \to D[1]$, 
then 
\[
m_{\sigma,t}(E) \le m_{\sigma,t}(D)+m_{\sigma,t}(F).
\]
\end{prop}
The proof of Proposition \ref{prop:mass_triangle}
 is given in Section \ref{sec:tri}. 
 
\begin{prop}
\label{prop:bound}
Let $\sigma$ be a stability condition on $\D$. 
Then 
\[
m_{\sigma,t}(E)\le m_{\sigma,t}(D)\delta_t(D,E)
\]
for any objects $D,E \in \D$. 
\end{prop}
{\bf Proof. } It is sufficient to show the case $E \in \left<D\right>$. 
Then by the definition of complexity $\delta_t(D,E)$,
for any $\eps>0$ there is a sequence of exact triangles
\[
\xymatrix @C=4.5mm{ 
 0 \ar[rr]   &&  A_1 \ar[dl] \ar[rr] && A_2 \ar[dl] 
 \ar[r] & \dots  \ar[r] & A_{k-1} \ar[rr] && E\oplus E^{\prime} \ar[dl] \\
& D[n_1] \ar@{-->}[ul] && D[n_2] \ar@{-->}[ul] &&&& D[n_k] \ar@{-->}[ul] 
}
\]
such that 
\[
\sum_{i=1}^k e^{n_i t}<\delta_t(D,E)+\eps.
\]
Note that $m_{\sigma,t}$ satisfies $m_{\sigma,t}(D[n])=
m_{\sigma,t}(D)\cdot e^{nt}$ for $ D \in \D$ and $n \in \Z$.
By using the inequality in Proposition \ref{prop:mass_triangle} repeatedly, 
we have
\begin{align*}
m_{\sigma,t}(E) \le m_{\sigma,t}(E\oplus E^{\prime})
&\le \sum_{i=1}^k m_{\sigma,t}(D[n_i]) 
\le m_{\sigma,t}(D)\cdot \left(\sum_{i=1}^k e^{n_i t}\right)  \\
&\le m_{\sigma,t}(D)\delta_t(D,E)+\eps \cdot m_{\sigma,t}(D)
\end{align*}
for any $\eps>0$. This implies the result. 
\endproof

Now we show fundamental properties of  the mass growth.
\begin{thm}
\label{thm:well}
Let $F \colon \D \to \D$ be an endofunctor and $\sigma$ be 
a stability condition on $\D$. Assume that $\D$ has a 
split-generator $G \in \D$. 
Then the mass growth $h_{\sigma,t}(F)$ satisfies the followings. 
\begin{itemize}
\item[$(1)$] $h_{\sigma,t}(F)$ is given by
\[
h_{\sigma,t}(F)=\limsup_{n \to \infty}\frac{1}{n}\,
\log (m_{\sigma,t}(F^n G)).
\]
\item[$(2)$] We have an inequality 
\[
h_{\sigma,t}(F) \le h_t(F) < \infty
\]
where $h_t(F)$ is the entropy of $F$ (see Definition \ref{def:entropy}). 
\end{itemize}
\end{thm}
{\bf Proof. } By Proposition \ref{prop:DHKK} $(3)$ 
and Proposition \ref{prop:bound}, we have
\[
m_t(F^n E) \le m_t(F^n G)\delta_t(F^n G,F^n E) 
\le m_t(F^n G)\delta_t(G, E)
\] 
for any object $E \in \D$. 
Hence 
\[
\limsup_{n  \to \infty} \frac{1}{n}\,\log m_t(F^n E) \le
\limsup_{n  \to \infty} \frac{1}{n}\,\log m_t(F^n G)
\]
and this inequality implies $(1)$. 
Again by Proposition \ref{prop:bound}, we have
\[
m_t(F^n G) \le m_t(G) \delta_t(G,F^n G). 
\]
Hence
\[
\limsup_{n  \to \infty} \frac{1}{n}\,\log m_t(F^n G) \le
\lim_{n  \to \infty} \frac{1}{n}\,\log \delta_t(G,F^n G)
\]
and this inequality implies $(2)$. 
\endproof

\subsection{Triangle inequality for mass with a parameter}
\label{sec:tri}
We prove Proposition \ref{prop:mass_triangle}. 
Recall the notation $\bH=\{\,r e^{i \pi \phi}\,|\,r>0,\, \phi \in (0,1]\,\}$.
For a complex number $z \in \bH$, define the 
function of $t \in \R$
by 
\[
g_t(z):=|z|\,e^{\phi(z)t}
\]
where $\phi(z)$ is the phase of $z$ given by  
$\phi(z):=(1 \slash \pi)\arg z \in (0,1]$.
We start to show the triangle inequality for $g_t(z)$. 
\begin{lem}
For $z_1,z_2 \in \bH$, an inequality
\[
g_t(z_1+z_2) \le g_t(z_1)+g_t(z_2)
\]
holds. 
\end{lem}
\begin{figure}[!ht]
\centering
\includegraphics[scale=1.1]{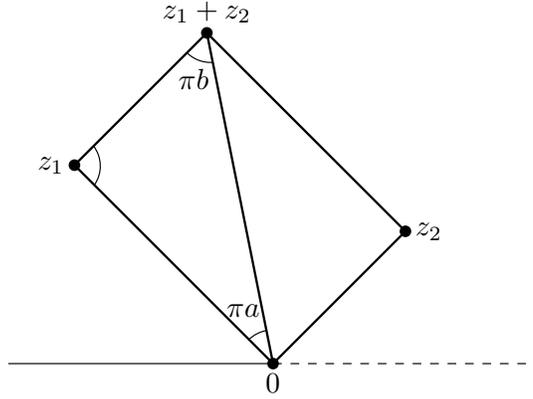}
\caption{Triangle consisting of vertices $0,z_1, z_1+z_2$.}
\label{Fig:Triangle}
\end{figure}
{\bf Proof.} 
Set $\phi_1:=\phi(z_1), \phi_2:=\phi(z_2)$ and $\phi_3:=\phi(z_1+z_2)$.  
If $\phi_1=\phi_2$, the result is trivial. We consider the case $\phi_1>\phi_2$. 
By applying the law of sine for the triangle consisting of 
vertices $0, z_1, z_1+z_2$ (see Figure \ref{Fig:Triangle}), we obtain
\[
|z_1+z_2|=d \sin (\pi a+\pi b),\quad
|z_1|=d \sin \pi b,\quad |z_2|=d \sin \pi a 
\]
where $a=\phi_1-\phi_3$, $b=\phi_3-\phi_2$ 
and $d$ is the diameter of the excircle of the triangle. 
By inputting these parameters, 
the inequality $g_t(z_1+z_2) \le g_t(z_1)+g_t(z_2)$ becomes
\[
\sin(\pi a+\pi b)\le e^{at}\sin \pi b+e^{-bt}\sin \pi a
\]
where $0<a,b<1$. Dividing by $\sin \pi a \sin \pi b$ and applying 
the addition formula, we have
\[
\frac{e^{at}-\cos \pi a  }{\sin \pi a}+\frac{e^{-bt}-\cos \pi b  }{\sin \pi b}
\ge 0.
\]
After setting $c=-b$, the above inequality is 
equivalent to $f(a) \ge f(c)$
for $-1<c<0<a<1$ where
\[
f(x)=\frac{e^{xt }-\cos \pi x  }{\sin \pi x}. 
\]
It is easy to check that $f(x)$ is  increasing in the intervals $(-1,0)$ 
and $(0,1)$, and the limit of  $f(x)$ at the zero 
is given by 
$\lim_{x \to \pm 0}f(x)=\frac{t}{\pi}$.
\endproof

\begin{figure}[!ht]
\centering
\includegraphics[scale=1]{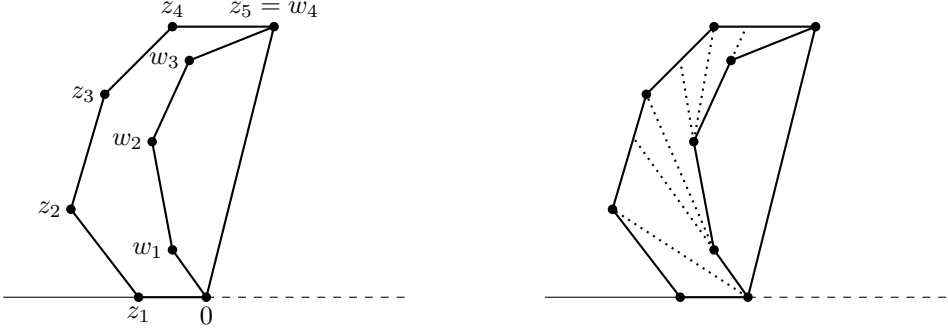}
\caption{Polygons and a triangulation of the encircled domain.}
\label{Fig:Polygons}
\end{figure}

The triangle inequality for $g_t(z)$ implies the following. 
\begin{lem}
\label{lem:polygon}
Let $z_1,\dots,z_k$ and $w_1,\dots,w_l$ be complex numbers in $\bH$ 
with $z_k=w_l$ and set $z_0=w_0=0$. 
If they satisfy the following conditions (see the left of Figure \ref{Fig:Polygons}):
\begin{itemize} 
\item[(a)] $\phi(z_i-z_{i-1})>\phi(z_{i+1}-z_i)$ and 
$\phi(w_j-w_{j-1})>\phi(w_{j+1}-w_j)$ 
for $i=1,\dots,k$ and $j=1,\dots,l$,
\item[(b)] the polygon $w_0w_1w_2 \dots w_l w_0$ contains the polygon 
$z_0 z_1 z_2 \dots z_k z_0$,
\end{itemize}
then
\[
\sum_{i=1}^k g_t(z_i-z_{i-1}) \le \sum_{j=1}^l g_t(w_j-w_{j-1}).
\]
\end{lem}
{\bf Proof. }
By the condition $(b)$, there is a unique domain encircled by two paths 
$z_0z_1z_2 \dots z_k$ and $w_0 w_1 w_2 \dots w_l$. By the convexity condition $(a)$, 
we can triangulate this domain as in the right of 
Figure \ref{Fig:Polygons}. Applying the triangle inequality for $g_t(z)$ (Lemma \ref{sec:tri}) 
repeatedly, we obtain the result.
\endproof

\begin{lem}
\label{lem:mass_triangle}
Let $\sigma=(Z,\sli)$ be a stability condition and $\H=\sli((0,1])$ be the 
associated heart. If there is a short exact sequence 
\[
0 \to A \to B \to C \to 0
\]
in $\H$ and $C \in \sli(1) $, 
then
\[
m_t(A) \le m_t(B)+e^{-t}m_t(C). 
\]
\end{lem}
{\bf Proof. } Let 
\begin{align*}
&0=A_0 \subset A_1 \subset A_2 \subset \cdots \subset A_l=A \\
&0=B_0 \subset B_1 \subset B_2 \subset \cdots \subset B_{k-1}=B
\end{align*}
be Harder-Narasimhan filtrations of $A$ and $B$.  
Set $z_i:=Z(A_i)$ for 
$i=0,1,\dots,k$, $w_j:=Z(B_j)$ for $j=0,1,\dots,l-1$ and $w_l:=Z(B)-Z(C)=Z(A)$. 
Then by definition of the Harder-Narasimhan filtration, these complex numbers
 satisfy the condition $(a)$ in Lemma \ref{lem:polygon}. Consider 
 the Harder-Narasimhan polygons $\HN^Z(A)$ and $\HN^Z(B)$ 
 (see Definition \ref{defi:HN_polygon}). By Proposition \ref{prop:HN_polygon}, 
 complex numbers $z_0,z_1,\dots,z_k$ and $w_0,w_1,\dots,w_{l-1}$ are 
 extremal points of $\HN^Z(A)$ and $\HN^Z(B)$ respectively. 
 Thus the intersection of $\HN^Z(A)$ and the left of the 
 line through $0$ and $z_k=Z(A)$ is the polygon $z_0z_1z_2 \dots z_k z_0$
 and the intersection of $\HN^Z(B)$ and the left of the 
 line through $0$ and $w_l=Z(A)$ is the polygon $w_0w_1w_2 \dots w_l w_0$. 
 Since $\HN^Z(A) \subset \HN^Z(B)$, the polygon $w_0w_1w_2 \dots w_l w_0$
contains the polygon $z_0z_1z_2 \dots z_k z_0$  and this implies the 
condition $(b)$ in Lemma \ref{lem:polygon}.
Since
\begin{align*}
m_t(A)=\sum_{i=1}^{k}g_t(z_i-z_{i-1}), \quad
m_t(B)=\sum_{j=1}^{l-1} g_t(w_j-w_{j-1}), \quad  
e^{-t} m_t(C)=g_t(w_l-w_{l-1}),
\end{align*}
applying Lemma \ref{lem:polygon}, we obtain the result. 
\endproof

{\bf Proof of Proposition \ref{prop:mass_triangle}. }
Assume that there is a exact triangle $D \to E \to F\to D[1]$. 
From a Harder-Narasimhan filtration of $E$, we can construct 
the dual Harder-Narasimhan filtration
\[
D =
\xymatrix @C=4.5mm{ 
 D_m \ar[rr]   &&  D_{m-1} \ar@{-->}[dl] \ar[rr] && D_{m-2} \ar@{-->}[dl] 
 \ar[r] & \dots  \ar[r] & D_{1} \ar[rr] && D_0 \ar@{-->}[dl] \\
& A_m \ar[ul] && A_{m-1} \ar[ul] &&&& A_1 \ar[ul] 
}
=0
\]
with $A_i \in \sli(\phi_i)$ and $\phi_m>\phi_{m-1}>\cdots \phi_1$. Applying 
the octahedra axiom for the above sequence together with 
the exact triangle $D \to E \to F\to D[1]$, we can construct a sequence of 
exact triangles
\[
E=
\xymatrix @C=4.5mm{ 
E_m \ar[rr]   &&  E_{m-1} \ar@{-->}[dl] \ar[rr] && E_{m-2} \ar@{-->}[dl] 
 \ar[r] & \dots  \ar[r] & E_{1} \ar[rr] && F \ar@{-->}[dl] \\
& A_m \ar[ul] && A_{m-1} \ar[ul] &&&& A_1 \ar[ul]
}.
\]
Since $m_t(D)=\sum_{i=1}^m |Z(A_i)|e^{t \phi_i}$ and $A_i$ is semistable, 
the problem is reduced to the case that $D$ is semistable. 
Without loss of generality we can assume $D \in \sli(1)$. By taking 
the cohomology associated with the heart $\H=\sli((0,1])$ 
(see Section \ref{sec:t-structure}), we have 
a long exact sequence 
\[
0 \to H^{-1}(E) \to H^{-1}(F) \to H^0(D) \to H^0(E) \to H^0(F) \to 0
\]
and isomorphisms $H^i(E) \cong H^i(F)$ for $i \neq -1,0$ in $\H$.  
If $1>\phi^+(H^0(E))$, then the map $f\colon H^0(D) \to H^0(E)$ is zero. Hence 
the long exact sequence splits into 
$0 \to H^{-1}(E) \to H^{-1}(F) \to H^0(D) \to 0$ and $H^0(E) \cong H^0(F)$. 
From Lemma \ref{lem:mass_triangle}, we have 
\[
m_t(H^{-1}(E))e^t \le  m_t(H^{-1}(F)) e^t +m_t(D).
\]
Thus we obtain the result. If  the maps $f \colon H^0(D) \to H^0(E)$
is not zero, then the long exact sequence splits into two short exact sequences
\begin{align*}
0 \to H^{-1}(E) \to &H^{-1}(F) \to \Ker f \to 0  \\ 
 0 \to\Im f \to &H^0(E) \to H^0(F) \to 0.
\end{align*}
Let $E_+ \in \sli(1)$ the semistable factor of $E$ with the phase one. 
Note that $m_t(D)=m_t(\Ker f)+m_t(\Im f) $
 since $\Ker f \subset D \in \sli(1)$ and  $\Im f \subset E_+ \in \sli(1)$.   
Again by Lemma \ref{lem:mass_triangle}, 
we have 
\[
m_t(H^{-1}(E))e^t \le  m_t(H^{-1}(F)) e^t +m_t(\Ker f)
\]
and it is easy to check that $m_t(H^0(E))=m_t(\Im f)+m_t(H^0(F))$. 
\endproof

\subsection{Mass growth and deformation of stability conditions}
The aim of this section is to show that for a stability condition $\sigma$ 
and an endofunctor $F$, the mass growth
$h_{\sigma,t}(F)$ is stable under the continuous deformation of $\sigma$. 
The following inequality was shown in \cite[Proposition 8.1]{Br1}
when $t=0$. 
\begin{prop}
\label{prop:mass_inequality}
Let $\sigma=(Z,\sli) \in \Stab(\D)$ be a stability condition on $\D$. 
If $\tau=(W,\Q) \in B_{\eps}(\sigma)$ with small enough $\eps>0$, 
then there are functions $C_1,C_2 \colon \R \to \R_{>0}$ such that
\[
C_1(t) \,m_{\tau,t}(E)<m_{\sigma,t}(E)<C_2(t)\,m_{\tau,t}(E)
\]
for all $0 \neq E \in \D$.
\end{prop}
{\bf Proof. }We use the argument similar to  the proof of 
\cite[Proposition 8.1]{Br1}. 
It is sufficient to show that for $\tau=(W,\Q) \in B_{\eps}(\sigma)$ 
with small enough $\eps>0$, there is some constants $C>1$ and $r>0$ 
such that 
\[
m_{\tau,t}(E)<C e^{r |t|}\, m_{\sigma,t}(E)
\] 
for any nonzero object $E \in \D$. 
We first consider the case $\phi_{\sigma}^+(E)-\phi_{\sigma}^-(E)<\eta$ 
for $0<\eta<\frac{1}{4}$. 
In this case, it was shown in the proof of \cite[Proposition 8.1]{Br1}
that there is a constant $C(\eps,\eta)>1$ such that 
\[
m_{\tau}(E) \le C(\eps,\eta)\,m_{\sigma}(E)
\]
and $C(\eps,\eta)\to 1$ as $\max\{\eps,\eta\} \to 0$. 
Note that $\phi_{\sigma}^+(E)-\phi_{\sigma}^-(E)<\eta$ implies 
$\phi_{\sigma}^{\pm}(E) \in (\psi,\psi+\eta)$ for some $\psi \in \R$.
Since $d(\sli,\Q)<\eps$, we have 
$\phi_{\tau}^{\pm}(E) \in (\psi-\eps,\psi+\eps+\eta)$. 
By definition of $m_{\sigma,t}(E)$ and $m_{\tau,t}(E)$, it follows that
\[
m_{\tau,t}(E)\le  m_{\tau}(E)\exp\left({\phi_{\tau}^+(E)}|t|\right),\quad
m_{\sigma}(E)\exp\left({\phi_{\sigma}^-(E)}|t|\right)\le  m_{\sigma,t}(E).
\]
Since $\psi<\phi_{\sigma}^-(E)$ and $\phi_{\tau}^+(E)<\psi+\eps+\eta$, 
we have an inequality
\[
m_{\tau,t}(E) \le C(\eps,\eta) e^{(\eps+\eta)|t|} \,m_{\sigma,t}(E).
\]
Next we consider a general nonzero object $E$. Take a real number 
$\phi$ and a positive integer $n$. For $k \in \Z$, define intervals
\[
I_k:=[\phi+kn\eps,\phi+(k+1)n\eps),\quad
J_k:=[\phi+kn\eps-\eps,\phi+(k+1)n\eps+\eps)
\]
and let $\alp_k$ and $\beta_k$ be the truncation functors projecting into 
the subcategories $\Q(I_k)$ and $\sli(J_k)$ respectively. 
Again by the argument in the proof of \cite[Proposition 8.1]{Br1}, 
for  small enough $n\eps$,
we have
\[
m_{\tau,t}(E)=\sum_k m_{\tau,t}(\alp_k(E))
\le \sum_k m_{\tau,t}(\beta_k(E))<C(\eps,(n+2)\eps)e^{(n+3)\eps|t|}
\sum_k m_{\sigma,t}(\beta_k(E)). 
\]
On the other hand, we can choose $\phi$ so that
\[
\sum_k m_{\sigma,t}(\beta_k(E)) \le \frac{n+2}{n}m_{\sigma,t}(E). 
\]
By taking the limits $\eps \to 0$ and $n \to \infty$ 
in keeping with $n \eps \to 0$, 
the result follows. 
\endproof

From Proposition \ref{prop:mass_inequality}, 
we immediately have the following. 
\begin{prop}
\label{prop:connected}
Let $F \colon \D \to \D$ be an endofunctor, and 
$\sigma $ and $\tau$ be stability conditions on $\D$.  
If $\sigma$ and $\tau$ lie in the same connected component in $\Stab(\D)$, 
then 
\[
h_{\sigma,t}(F)=h_{\tau,t}(F).
\]
\end{prop}
{\bf Proof. } Let $\sigma,\tau \in \Stab(\D)$ be stability conditions 
such that $\tau \in B_{\eps}(\sigma)$ for small enough $\eps>0$. 
Then Proposition \ref{prop:mass_inequality} implies 
$h_{\sigma,t}(F)=h_{\tau,t}(F)$. Thus $h_{\sigma,t}(F)$ is 
locally constant on the topological space $\Stab(\D)$. 
\endproof

\subsection{Lower bound of the mass growth by the spectral radius}
Let $F \colon \D \to \D$ be an endofunctor. 
Since $F$ preserves exact triangles in $\D$, 
$F$ induces a linear map 
\[
[F]\colon K(\D) \to K(\D).
\]
The {\it spectral radius of $[F]$} is defined by
\[
\rho([F]):=\max \{|\lam|\,|\, \lam \text{ is an eigenvalue of }[F]  \,   \}. 
\]
\begin{prop}
\label{prop:spectral}
For any stability condition $\sigma \in \Stab(\D)$, 
we have an inequality
\[
\log \rho([F]) \le h_{\sigma,0}(F).
\]
\end{prop}
{\bf Proof. } Set $K(\D)_{\C}:=K(\D) \otimes \C$. 
Let $A_1,\dots,A_n \in \D$ be objects whose classes 
$[A_1],\dots,[A_n]$ form a basis of $K(\D)_{\C}$. 
Take an eigenvector 
\[
v=\sum_{i=1}^n  a_i [A_i] \in K(\D)_{\C} \quad (\,a_i \in \C\,)
\] for the eigenvalue $\lam \in \C$ of $[F]$
 satisfying $|\lam|=\rho([F])$. 
First we consider the case that a stability condition $\sigma=(Z,\sli)$ 
satisfies $Z(v)\neq 0$. Note that the mass satisfies 
$ |Z(E)| \le  m_{\sigma}(E) $ and $m_{\sigma}(E \oplus E^{\prime})=
m_{\sigma}(E)+m_{\sigma}(E^{\prime})$ for any objects 
$E,E^{\prime} \in \D$. 
Then
\begin{align*}
|\lam|^k |Z(v)|=|Z(\lam^k v)|=|Z([F]^k v)| &\le
\sum_{i=1}^n |a_i|\cdot|Z(F^k A_i)|  \\
&\le \sum_{i=1}^n l_i\, m_{\sigma} (F^k A_i) 
=m_{\sigma}\left( F^k \left(\oplus_{i=1}^n A_i^{\oplus l_i}\right)  \right)
\end{align*}
where $l_1,\dots,l_n$ are positive integers satisfying $|a_i| \le l_i$. 
Since $|Z(v)|>0$, we have 
\[
\log\rho([F])=\lim_{k \to \infty}\frac{1}{k}\log(|\lam|^k|Z(v)|)
\le \limsup_{k \to \infty}\frac{1}{k} \log (m_{\sigma}(F^k E))
\le h_{\sigma,0}(F)
\]
where  $E=\oplus_{i=1}^n A_i^{\oplus l_i}$. 
Next consider the case $Z(v)=0$. Then by Theorem \ref{thm:localiso}, 
we can deform $\sigma=(Z,\sli)$ to $\sigma^{\prime}=(Z^{\prime},\sli^{\prime})$ 
so that $Z^{\prime}(v) \neq 0$. Again we have 
$\log \rho([F]) \le h_{\sigma^{\prime},0}(F)$ and 
Proposition \ref{prop:connected} implies 
$h_{\sigma,0}(F)=h_{\sigma^{\prime},0}(F)$.
\endproof

\subsection{Mass growth via algebraic stability conditions}
\label{sec:inequality_algebraic}
If a triangulated category has an algebraic stability condition,  
then we can show that the mass growth coincides with the entropy. 
Let $\H \subset \D$ be an algebraic heart with simple objects 
$S_1,\dots,S_n$. Then the Grothendieck group is given by 
\[
K(\D) \cong \bigoplus_{i=1}^n \Z [S_i]. 
\]
The class of an object $E \in \H$ is written as 
$[E]=\sum_{i=1}^n d_i [S_i]$
with $d_i \in \Z_{\ge 0}$. We define the {\it dimension of $E$} by 
$\dim E:=\sum_{i=1}^n d_i \in \Z_{\ge 0}$. Then the dimension gives the upper 
bound of the complexity for objects in $\H$.  
\begin{lem}
\label{lem:dim}
Let $\H \subset \D$ be an algebraic heart with simple objects 
$S_1,\dots,S_n$. Then for the split-generator $G:=\oplus_{i=1}^n S_i$,  
we have an inequality
\[
\delta_{t}(G,E) \le \dim E.
\]
\end{lem}
{\bf Proof. } Since $\H$ is a finite length abelian category, 
for any object $E \in \H$ there is a Jordan-H\"older 
filtration 
\[
0=E_0 \subset E_1 \subset E_2 \subset \cdots  \subset E_l=E
\]
of length $l=\dim E$  with $E_i \slash E_{i-1} \in \{S_1,\dots,S_n\}$. 
As a result, we can construct a filtration
\[
0=E_0^{\prime} \subset E_1^{\prime} \subset E_2^{\prime} 
\subset \cdots  \subset E_l^{\prime}=E \oplus E^{\prime}
\]
of length $l=\dim E$ with $E_{i}^{\prime}\slash E_{i-1}^{\prime}=G$
and this implies the result.
\endproof

Following Section \ref{sec:algebraic_stab}, 
we construct the special algebraic stability condition. 
For an algebraic heart $\H \subset \D$ with simple objects 
$S_1,\dots,S_n$, define the central charge 
\[
Z_0 \colon K(\D)\cong \bigoplus_{i=1}^n \Z [S_i] \to \C
\]
by $Z_0(S_i):=i$. Then the pair $\sigma_0:=(Z_0,\H)$ 
becomes an algebraic stability condition. By definition,  
the mass of an object $E \in \H$ is given by 
\[
m_{\sigma_0,t}(E)=\dim E \cdot e^{\frac{1}{2}t}. 
\]
Together with Lemma \ref{lem:dim}, we obtain the following inequality. 
\begin{prop}
\label{prop:algebraic_bound}
For the generator $G=\oplus_{i=1}^n S_i$ and 
the algebraic stability condition  $\sigma_0=(Z_0,\H)$, we have an inequality
\[
\delta_t(G,E) \le e^{-\frac{1}{2}t} \, m_{\sigma_0,t}(E).
\]
\end{prop}
{\bf Proof. }
For an object $E \in \D$, we denote by $H^k(E) \in \H$ the cohomology 
associated with the heart $\H$ 
(see Section \ref{sec:t-structure}).
By using Lemma \ref{lem:triangle} and Lemma \ref{lem:dim}, we have
\begin{align*}
\delta_t(G,E)  \le \sum_{k} \delta_t(G,H^k(E)  )\,e^{-k t} 
\le \sum_{k} \dim H^k(E)  \,e^{-k t}.
\end{align*}
On the other hand, the definition of  $m_{\sigma_0,t}(E)$ 
implies
\[
m_{\sigma_0,t}(E)=\sum_k  m_{\sigma_0,t}(H^k(E))\,e^{-kt}
=\sum_k  \dim H^k(E) \,e^{\frac{1}{2}t}  \,e^{-kt}.
\] 
Thus we obtain the result. 
\endproof

We show the main result of this section. 
\begin{thm}
\label{thm:algebraic}
Let $G \in \D$ be a split-generator and $F \colon \D \to \D$ be 
an endofunctor. 
If a connected component 
$\Stab^{\circ}(\D) \subset \Stab(\D)$ contains an algebraic
stability condition, then for any $\sigma \in \Stab^{\circ}(\D)$
we have
\[
h_t(F)=h_{\sigma,t}(F)=\lim_{n \to \infty}\frac{1}{n}\log (m_{\sigma,t}(F^n G)).
\]
\end{thm}
{\bf Proof.}
Let $\H$ be an algebraic heart with simple objects $S_1,\dots,S_n$ and 
set $G=\oplus_{i=1}^n S_i$. Consider the special 
algebraic stability condition $\sigma_0=(Z_0,\H)$ 
which is constructed in this section.
By Proposition \ref{prop:connected}, it is sufficient to show that
\[
h_{\sigma_0,t}(F)=h_t(F). 
\]
By \cite[Lemma 2.6]{DHKK}, the limit
\[
h_t(F)=\lim_{n \to \infty}\frac{1}{n}\,\log \delta_t(G,F^n G)
\]
converges. On the other hand, by Proposition \ref{prop:bound}  
and Proposition \ref{prop:algebraic_bound}, we have
\[
e^{\frac{1}{2}t}\, \delta_t(G,F^n G)  \le
m_{\sigma_0,t}(F^n G) \le m_{\sigma_0,t}(G)\delta_t(G,F^n G).
\] 
Hence the limit
\[
\lim_{n \to \infty}\frac{1}{n}\,\log (m_{\sigma_0,t}(F^n G))
\]
converges and coincides with $h_t(F)$.
\endproof

\section{Applications}
\subsection{Entropy on the derived categories of non-positive dg-algebras}
\label{sec:dga}
In this section, we discuss the entropy of endofunctors 
on the derived categories of non-positive dg-algebras. In this case, 
we can describe the entropy as the growth rate of the 
Hilbert-Poincar\'e polynomial of a generator. 

Let $A =\oplus_{k \in \Z}A^k$ be a dg-algebra over $\K$ satisfying 
the following conditions:
\begin{itemize}
\item[(a)] $H^k(A)=0$ for $i>0$,
\item[(b)] $H^0(A)$ is a finite dimensional algebra over $\K$.
\end{itemize}

Let $\D(A)$ be the derived category of dg-modules over $A$ and 
 $\D_{fd}(A)$ be the full subcategory of $\D(A)$ consisting of dg-modules 
with finite dimensional total cohomology, i.e. 
\[
\D_{fd}(A):=\left\{\,M \in \D(A)\,\middle\vert\, 
\sum_k \dim H^k(M)<\infty \,\right\}.
\]
Define the full subcategory $\F \subset \D_{fd}(A)$ by
\[
\F:=\{\,  M \in \D_{fd}(A)\,\vert\,  H^k(M)=0 \text{ for }k>0 \, \}.
\] 
Then $\F$ becomes a bounded t-structure on $\D_{fd}(A)$.
The heart $\H_s$ of $\F$ is called the {\it standard heart}. 
It is known that the $0$-th cohomology functor 
 $H^0 \colon \D_{fd}(A) \to \Mod H^0(A)$ 
induces an equivalence of abelian categories:
\[
 H^0 \colon \H_s \iso \Mod H^0(A)
\]
where $\Mod H^0(A)$ is an abelian category of finite dimensional $H^0(A)$-modules.
(For details, see \cite[Section 2]{Am}.) Since $H^0(A)$ is a finite dimensional algebra, 
$\H_s$ becomes an algebraic heart. 
Thus we can construct an algebraic stability condition on $\D_{fd}(A)$. 
Applying Theorem \ref{thm:algebraic}, we obtain the following. 
\begin{prop}
Let $\Stab^{\circ}(\D_{fd}(A))$ be the connected component which contains 
stability conditions with the standard heart $\H_s$. Then for any stability conditions 
$\sigma \in \Stab^{\circ}(\D_{fd}(A))$ and an endofunctor $F \colon 
\D_{fd}(A)\to \D_{fd}(A)$, we have
\[
h_t(F)=h_{\sigma,t}(F). 
\]
\end{prop}
Next we describe $h_t(F)$ by using the Hillbert-Poincar\'e polynomial.
\begin{defi}
For a dg-module $M \in \D_{fd}(A)$, define 
the {\it Hilbert-Poincar\'e polynomial of $M$} 
by 
\[
P_t(M):=\sum_{k \in \Z}  \dim H^k(M) \,e^{-kt} \in \Z[e^{t},e^{-t}]. 
\]  
\end{defi}
As in Section \ref{sec:inequality_algebraic}, we construct the special stability 
condition $\sigma_0=(Z_0,\H_s)$ by using the standard heart $\H_s$. 
Then by definition of $\sigma_0$, we have
\[
m_{\sigma_0,t}(M)=e^{\frac{1}{2}t} P_t(M)
\]
for any dg-module $M \in \D_{fd}(A)$. 
As a result, the entropy is described as follows. 
\begin{prop}
\label{prop:HP}
Let $F \colon \D_{fd}(A) \to \D_{fd}(A)$ be an endofunctor and 
$G \in \D_{fd}(A)$ be a split-generator. Then the entropy of $F$ 
is given by
\[
h_t(F)=\lim_{n \to \infty}\frac{1}{n}\,\log P_t(F^n G).
\]
\end{prop}

\subsection{Entropy of spherical twists}
\label{sec:spherical}
In this section, we compute the entropy of 
Seidel-Thomas spherical twists on the derived categories of 
Calabi-Yau algebras associated with acyclic quivers. 
Let $Q$ be an acyclic quiver with vertices $\{1,\dots,n\}$
 and $\Gamma_N Q$ be the {\it 
Ginzburg $N$-Calabi-Yau dg-algebra associated with $Q$} for $N \ge 2$.  
(For the definition of $\Gamma_N Q$, 
see \cite[Section 4.2]{Ginz} or \cite[Section 6.2]{Kel}.) 
Set $\D_Q^N:=\D_{fd}(\Gamma_N Q)$.  
By \cite[Theorem 6.3]{Kel}, the category $\D_Q^N$ 
becomes a {\it $N$-Calabi-Yau category}, i.e. 
there is a natural isomorphism 
\[
\Hom(E,F) \iso \Hom(F,E[N])^*
\]
for $E,F\in \D_Q^N$. (The notation $V^*$ is the dual space of 
a $\K$-linear space $V$.)  
In the Calabi-Yau category, we can consider a certain 
class of objects, called spherical objects. 
An object $S \in \D_Q^N$ is called {\it $N$-spherical} if 
\begin{align*}
\Hom(S,S[i]) = 
\begin{cases}
\K \quad \text{if} \quad i=0,N  \\
\,0 \quad \text{otherwise} . 
\end{cases}
\end{align*} 
For a spherical object $S \in \D_Q^N$, Seidel-Thomas \cite{ST} 
defined an exact autoequivalence $\Phi_S \in \Aut(\D_Q^N)$, 
called a {\it spherical twist},  
by the exact triangle 
\[
\Hom^{\bullet}(S,E) \otimes S \lto E \lto \Phi_S(E)
\]
for any object $E \in \D^N_Q$.
The inverse functor $\Phi_S^{-1} \in \Aut(\D^N_Q)$ is given by
\[
\Phi_S^{-1}(E) \lto E \lto S \otimes \Hom^{\bullet}(E,S)^* .
\]
The Ginzburg dg-algebra $\Gamma_N Q$ satisfies the conditions in 
Section \ref{sec:dga} when $N \ge 2$. (In the case $N=2$, we need some modification.) 
Hence the category $\D_Q^N$ has the standard algebraic heart $\H_s \subset \D_Q^N$ 
generated by simple $\Gamma_N Q$-modules $S_1,\dots,S_n$ corresponding 
to vertices $\{1,\dots,n\}$ of $Q$.  In addition, these objects 
$S_1,\dots,S_n$ become $N$-spherical by \cite[Lemma 4.4]{Kel}. 
Thus we can define spherical twists  
$\Phi_{S_1},\dots,\Phi_{S_n} \in \Aut(\D_Q^N)$. 
In the following, we compute the entropy of 
spherical twists $\Phi_{S_1},\dots,\Phi_{S_n}$ 
by using  Proposition \ref{prop:HP}. 
For simplicity, write $\Phi_i:=\Phi_{S_i}$. 
\begin{lem}
\label{lem:computation}
For a spherical twist $\Phi_i\in \Aut(\D_Q^N)$ 
and a spherical object $S_j \in \D_Q^N$ , the Hilbert-Poincar\'e 
polynomial of $\Phi_i^k S_j\,(k \ge 0)$ 
is given by
\[
P_t(\Phi_i^k S_j) =
\begin{cases}
e^{k(1-N)t} &\text{if }i=j \\
1+q_{ij} \sum_{l=0}^{k-1}e^{l(1-N)t}  &\text{if }q_{ij}>0 \\
1+q_{ji} \,e^{(2-N)t} \sum_{l=0}^{k-1}e^{l(1-N)t} &\text{if }q_{ji}>0 \\
 1 &\text{otherwise }
\end{cases}
\]
where $q_{ij}$ is the number of arrows from $i$ to $j$ in $Q$.
\end{lem}
{\bf Proof.}
First we note that
\[
\dim \Hom(S_i,S_j[m])=
\begin{cases}
1 &\text{if }i=j \text{ and }m=0,N \\
q_{ij}&\text{if }q_{ij}>0 \text{ and }m=1 \\
q_{ji} &\text{if }q_{ji}>0 \text{ and }m=N-1 \\
 0 &\text{otherwise }.
\end{cases}
\]
By the definition of spherical twists, it is easy to see that 
$\Phi_i^k S_i=S_i[k(1-N)]$ and hence $P_t(\Phi_i^k S_i)=e^{k(1-N)t}$. 
If $i\neq j$ and $q_{ij}=q_{ji}=0$, then $\Phi^k_i S_j=S_j$ and hence 
$P_t(\Phi_i^k S_j)=1$. Consider the case $q_{ij}>0$. 
Since 
\[
\Hom^{\bullet}(S_i,S_j) \otimes S_i =\bigoplus_{m \in \Z}\Hom(S_i[m],S_j) 
\otimes S_i[m]
\cong S_i^{\oplus q_{ij}}[-1],
\]
we have an exact triangle
\[
S_j \to \Phi_i S_j \to S_i^{\oplus q_{ij}} \to S_j[1].
\]
Applying the spherical twist $\Phi_i$ for the above exact triangle repeatedly, 
we obtain a sequence of exact triangles
\[
\xymatrix @C=3mm{ 
 S_j \ar[rr]   &&  \Phi_i S_j \ar[dl] \ar[rr] && \Phi_i^2 S_j \ar[dl] 
 \ar[r] & \dots  \ar[r] & \Phi_i^{k-1}S_j \ar[rr] && \Phi_i^k S_j. \ar[dl] \\
& S_i^{\oplus q_{ij}} \ar@{-->}[ul] && S_i^{\oplus q_{ij}}[1-N] 
\ar@{-->}[ul] &&&& S_i^{\oplus q_{ij}}[(k-1)(1-N)] \ar@{-->}[ul]
}
\]
This implies the result in the case $q_{ij}>0$ 
and the similar argument gives the result in the case $q_{ji}>0$. 
\endproof

\begin{prop}
Let $Q$ be a connected acyclic quiver and assume that $Q$ is 
not a quiver with one vertex and no arrows. Then 
the entropy of spherical twists 
$\Phi_1,\dots,\Phi_n $ is given by 
\[
h_t(\Phi_i)=
\begin{cases}
0 &\text{if }t \ge 0 \\
(1-N)t &\text{if }t < 0.
\end{cases}
\]
\end{prop}
{\bf Proof. } 
We use the generator $G=\oplus_{j=1}^nS_j$. 
Then $P_t(\Phi_i^k G)=\sum_{j=1}^n P_t(\Phi_i^k S_j)$. 
Recall from Proposition \ref{lem:computation} that 
\[
P_t(\Phi_i^kS_j)=1+q_{ij} \sum_{l=0}^{k-1}e^{l(1-N)t}=1+q_{ij} 
\frac{1-e^{k(1-N)t}}{1-e^{(1-N)t}}
\]
in the case $q_{ij}>0$ and 
\[
P_t(\Phi_i^k S_j)=1+q_{ji} \,e^{(2-N)t} \sum_{l=0}^{k-1}e^{l(1-N)t}
=1+q_{ij}e^{(2-N)t}  \frac{1-e^{k(1-N)t}}{1-e^{(1-N)t}}
\]
in the case $q_{ji}>0$. First we consider the case $t>0$. 
Then the above two terms converge 
to some positive real numbers 
as $k \to \infty$ since $(1-N)t<0$. 
By the assumption of $Q$, the sum $\sum_{j=1}^n P_t(\Phi_i^k S_j)$ 
contains at least one of the above two. 
As a result, 
$\sum_{j=1}^n P_t(\Phi_i^k S_j)$ also converges 
to some positive real number as $k \to \infty$. Thus 
\[
h_t(\Phi_i)=\lim_{k \to \infty}\frac{1}{k} \log P_t(\Phi_i^k G)=0
\]
when $t>0$. 
Next consider the case $t<0$. Similarly we can show that 
$e^{-k(1-N)t}\sum_{j=1}^n P_t(\Phi_i^k S_j)$ 
converges to some positive real number as $k\to \infty$ since $-(1-N)t<0$. 
Thus
\begin{align*}
h_t(\Phi_i)&=\lim_{k \to \infty}\frac{1}{k} \log P_t(\Phi_i^k G)
=\lim_{k \to \infty}\frac{1}{k} \log e^{k(1-N)t}  e^{-k(1-N)t} P_t(\Phi_i^k G) \\
&=(1-N)t+ \lim_{k \to \infty}\frac{1}{k} \log   e^{-k(1-N)t} P_t(\Phi_i^k G)
=(1-N)t
\end{align*}
when $t<0$.
Finally we can easily check that $h_t(\Phi_i)=0$ when $t=0$. 
\endproof

\begin{rem}
The subgroup of autoequivalences generated by spherical twists
\[
\Sph(\D_Q^N):=\left<\Phi_1,\dots,\Phi_n\right>
\subset  \Aut(\D_Q^N)
\]
is called the Seidel-Thomas braid group. 
Here we only computed the 
entropy of generators $\Phi_1,\dots,\Phi_n$. It is important problem 
to compute the entropy $h_t(\Phi)$ for a general element $\Phi \in \Sph(\D_Q^N)$. 
\end{rem}

\subsection{Lower bound of the entropy on the derived categories of 
surfaces}
Let $X$ be a smooth projective variety over $\C$ and 
denote by $\Db(X)$ the bounded derived 
category of coherent sheaves on $X$. 
Define the Euler form 
$\chi \colon K(\Db(X)) \times K(\Db(X)) \to \Z$
by 
\[
\chi(E,F):=\sum_{i \in \Z}(-1)^i \dim_{\C} \Hom_{\Db(X)}(E,F[i]).
\]
The {\it numerical Grothendieck group} $N(X)$ is the quotient 
of $K(\Db(X))$ by the radical of the Euler form $\chi$.  
Let $\End^{FM}(\Db(X))$ be the semi-group consisting 
of Fourier-Mukai type endofunctors. 
Since these endofunctors preserve the radical of $\chi$, they induce 
linear maps on $N(X)$, i.e. the semi-group homomorphism 
\[
\End^{FM}(\Db(X)) \to \End(N(X)),\quad F \mapsto [F]
\]
is well-defined (see \cite[Section 5.1]{KT}). A stability condition 
$\sigma=(Z,\sli)$ is called {\it numerical} if 
$Z \colon  K(\Db(X)) \to \C$ factors through 
the numerical Grothendieck group $N(X)$. 

In \cite{Br2,AB}, a numerical stability condition on $\Db(X)$ was constructed 
when $\dim_{\C} X=2$. Applying Theorem \ref{thm:well} and 
Proposition \ref{prop:spectral}, we obtain the following lower bound of the 
entropy. 
\begin{prop}
Let $X$ be a smooth projective surface over $\C$ 
and $F \colon \Db(X) \to \Db(X)$ 
be a Fourier-Mukai type endofunctor. Then 
\[
\log \rho ([F]) \le h_0(F)
\]
where $\rho([F]) $ is the spectral radius of the 
induced linear map $[F]\colon N(X) \to N(X)$ and $h_0(F)$ 
is the entropy of $F$ at $t=0$.
\end{prop}

\begin{rem}
Let $X$ be a smooth projective variety over $\C$.
In \cite{KT}, they conjectured that  
the equality $\log \rho([F])=h_0(F)$ holds   
for any autoequivalence $F \in \Aut(\Db(X))$.  
This conjecture was shown for a curve in \cite{Kik} and 
for a variety with the ample canonical bundle or the ample anti-canonical 
bundle in \cite{KT}. 
\end{rem}

\end{document}